\documentclass[12pt]{article}
\usepackage{fontenc}
\usepackage{amsmath}
\usepackage{amsfonts}
\usepackage{delarray}
\usepackage{theorem}
\newtheorem{Thm}{Theorem}
\newtheorem{Rq}{Remark}

\newtheorem{Pro}{Proposition}

\def\C{\mathbb C}
\def\Z{\mathbb Z}

\def\P{\mathbb P}
\newenvironment{proof}[1][Proof]{\begin{trivlist}
\item[\hskip \labelsep {\bfseries #1}]}{\end{trivlist}}

\title{Families of Painlev\'e VI  equations having a common solution}

\author{
Bassem \textsc{Ben Hamed} \footnote{ The first author was partially supported by a grant of the Faculty of
Sciences of Sfax,
Tunisia.} \\
Facult\'e des Sciences de Sfax, D\'epartement de
Math\'ematiques\\
BP 802, route Soukra km 4 Sfax 3018, Tunisia\\  \\
Lubomir \textsc{Gavrilov}\\
Universit\'e Paul Sabatier, MIG\\
Laboratoire Emile Picard,  UMR 5580\\
 31062 Toulouse Cedex 9, France }
\begin{document}
\maketitle
\begin{abstract}
We classify all functions  satisfying non-trivial families of $\mathbf{PVI_\alpha}$ equations. It turns out
that, up to an Okamoto equivalence, there are exactly four families parameterized by affine planes or lines.
Each affine space is generated by points of "geometric origin", associated either to deformations of elliptic
surfaces with four singular fibers, or to deformations of three-sheeted covers of $\P^1$ with branching locus
consisting of four points.
\end{abstract}
%Key words: isomonodromic deformation, Painlev\'e 6 equation, Picard-Fuchs equation.
\section{Introduction}

Consider the Painlev\'e VI ( $ \mathbf{PVI_\alpha}$) equation
\begin{eqnarray*}
 \frac{d^2 \lambda}{dt^2} &= &\frac{1}{2}(\frac{1}{\lambda}+ \frac{1}
{\lambda - 1} + \frac{1}{\lambda - t})(\frac{d\lambda}{dt})^2 - (\frac{1}{t} + \frac{1}{t - 1}
 + \frac{1}{\lambda - t})\frac{d\lambda}{dt}\\
& &+ \frac{\lambda (\lambda - 1)(\lambda - t)}
 {t^2 (t - 1)^2} [ \alpha_0 - \alpha_1 \frac{t}{\lambda ^2} + \alpha_2 \frac{t - 1}{(\lambda - 1)^2}
 +  (\frac{1}{2} - \alpha_3) \frac{t(t - 1)}{(\lambda - t)^2} ]  .
\end{eqnarray*}
parameterized by $\alpha = (\alpha_0,\alpha_1,\alpha_2,\alpha_3) \in \mathbb{C}^4$. Although any solution  of $
\mathbf{PVI_\alpha}$, for generic $\alpha_i$, is transcendental (and even a "new transcendental function") there
is a large amount of solutions which are algebraic in $t$. Their general classification is still an open problem
(e.g. \cite[Manin]{10}), except in the particular case $\alpha_1=\alpha_2=\alpha_3=0$ \cite[Dubrovin,
Mazzocco]{dub00}. The present paper addresses the question of classifying \emph{families} of algebraic
solutions. The simplest case occurs when a given algebraic solution satisfies each member of a non-trivial
family of $ \mathbf{PVI_\alpha}$ equations. By a non-trivial family of $ \mathbf{PVI_\alpha}$ equations we mean
a set $\{\mathbf{PVI_\alpha}\}_\alpha$ containing at least two distinct elements corresponding to, say,
$\alpha'$ and $\alpha''$. Then this solution satisfies the $ \mathbf{PVI_\alpha}$ equations corresponding to the
affine line containing $\alpha'$ and $\alpha''$. It follows that  each non-trivial family as above corresponds
to an affine plane in the parameter space $\mathbb{C}^4\{\alpha\}$. We classify all such affine spaces, together
with their associated algebraic  solutions (Theorem \ref{t7}). The proof of Theorem  \ref{t7} does not use the
notion of Picard-Fuchs equation. It turns out that the solutions which we obtain
 coincide surprisingly  with the solutions obtained earlier by Doran, who used deformations of elliptic
surfaces with four singular fibers and the related Picard-Fuchs equations, see Theorem \ref{t2}.

The second purpose of the paper is to give a (partial) explanation of the above coincidence.
 It is well known that each solution $(\lambda(t),\alpha)$ of a given $
\mathbf{PVI_\alpha}$ equation governs the isomonodromy deformation of a $2 \times 2$ Fuchsian system with four
singular points. We say that such a deformation is "geometric", if there is a fundamental matrix of solutions
whose entries are Abelian integrals depending algebraically on the deformation parameter, see
\cite[Doran]{doran}. A geometric deformation of a Fuchsian system is isomonodromic and the associated solution
of the corresponding Schlesinger system (or $ \mathbf{PVI_\alpha}$ equation) is an algebraic function in $t$.
When this holds true, we  say that the algebraic solution $(\lambda(t),\alpha)$ of $ \mathbf{PVI_\alpha}$ is of
geometric origin. The solutions found by Doran are of geometric origin and they correspond to special values of
$\alpha$ which belong to the affine spaces of  $ \mathbf{PVI_\alpha}$ equations described in Theorem \ref{t7}.
We shall prove that the same list of solutions can be obtained from deformations of ramified covers of $\P^1$
with four ramification points. The corresponding values of the parameters $\alpha$ are different and are shown
on Table \ref{table5}, see Theorem \ref{t3}. Finally we note that all these points of geometric origin generate
the affine planes of $ \mathbf{PVI_\alpha}$ equations, described in Theorem \ref{t7}.

\emph{Acknowledgments.}
 The solution
1A, as well its relation to the family 2A was pointed out to us by Philip Boalch. We acknowledge his interest
and critical remarks.

\section{Families of Painlev\'e VI  equations having a common solution}

Let $\lambda = \lambda(t)$ be a solution of the equations $\mathbf{PVI_{\alpha'}}$, $\mathbf{PVI_{\alpha''}}$,
$\alpha' \neq \alpha''$. Then $\lambda(t)$ satisfies the implicit equation
\begin{equation}\label{1}
\beta_0 - \beta_1 \frac{t}{\lambda ^2} + \beta_2 \frac{t - 1}{(\lambda - 1)^2}
  - \beta_3 \frac{t(t - 1)}{(\lambda - t)^2} = 0
\end{equation}
where $\beta =\alpha'-\alpha''=(\beta_0,\beta_1,\beta_2,\beta_3)$, and hence it is an algebraic function. The
function $\lambda(t)$ satisfies, moreover, the family  $\{ \mathbf{PVI_{\alpha}}\}_\alpha$, where $\alpha$
belongs to the affine line
$$
\{ \alpha' + s (\alpha' - \alpha''): s \in \mathbb{C} \} \subset \mathbb{C}^4 .
$$
It is seen from this that the set of all $\alpha$ such that $\mathbf{PVI_{\alpha}}$ is satisfied by the function
$\lambda(t)$ form an affine subspace of $\mathbb{C}^4$. We refer to the set of these $\mathbf{PVI_{\alpha}}$
equations as to a \emph{family of Painlev\'e VI  equations having a common solution}.

\begin{Thm}
\label{t7}
 The list of all families of Painlev\'e VI  equations having a common solution, together with the corresponding solution,
 is shown in Table~\ref{table1}.
\end{Thm}
 \begin{Rq}
Each solution $\lambda(t)$ is defined by a relation $P(\lambda(t),t)\equiv 0$ where $P$ is an irreducible
polynomial  given in Table \ref{table1}. The solutions in each of the 5 series of families on Table \ref{table1}
are equivalent up to a $S_4$-symmetry of Painlev\'e VI equation (see section \ref{s4}).

Moreover, the solutions  $1A,1B,\dots, 1F$ are Okamoto equivalent to the solutions $2A,2B,2C$. More precisely,
the solution $1A$
$$
a\lambda^2 - bt = 0, \alpha =(a,b,\frac{1}{8},\frac{1}{8})
$$
is equivalent, after applying the transformation $w_2$ of Okamoto \cite[p.363]{okamoto}, to the solution $
\lambda^2-t=0$ where the parameter $ \alpha $ equals to
 $$ ( (\frac{\sqrt{2a} -
\sqrt{2b}}{2\sqrt{2}})^2,(\frac{\sqrt{2a} - \sqrt{2b}}{2\sqrt{2}})^2,(1- \frac{\sqrt{2a} -
\sqrt{2b}}{2\sqrt{2}})^2,(1- \frac{\sqrt{2a} - \sqrt{2b}}{2\sqrt{2}})^2) .
$$
Thus, up to Okamoto equivalence, the families of Painlev\'e VI  equations having a common solution are
represented (for instance) by the four families $2A,3A,4A,5A$ on Table \ref{table1}.
\end{Rq}
\begin{Rq}
The functions $0, 1, t$ are not considered as solutions of Painlev\'e VI equation, and therefore are excluded
from Table \ref{table1}.
\end{Rq}
\begin{table}
\begin{center}
\begin{tabular}{||l|l|l||}
\hline\hline
Name & Solution of $\textbf{PVI}_{\alpha}$ equation & $\textbf{PVI}_{\alpha}$ equation  \\ \hline\hline
1A & $a\lambda^2-b t$ & $(a,b,\frac{1}{8},\frac{1}{8})$ \\ \hline
1B & $a(\lambda-1)^2+ b(t-1) $ & $(a,\frac{1}{8},b,\frac{1}{8})$ \\ \hline
1C & $a(\lambda-t)^2-b t(t-1)$ & $(a,\frac{1}{8},\frac{1}{8},b)$ \\ \hline
1D & $-at(\lambda-1)^2+ b (t-1)\lambda^2$ & $(\frac{1}{8},a,b,\frac{1}{8})$ \\ \hline
1E & $a(\lambda-t)^2+b (t-1) \lambda^2$ & $(\frac{1}{8},a,\frac{1}{8},b)$ \\ \hline
1F & $a(\lambda-t)^2-b t(\lambda-1)^2$ & $(\frac{1}{8},\frac{1}{8},a,b)$ \\ \hline \hline
2A & $\lambda^2 - t $ & $(a,a,b,b)$ \\ \hline
2B & $\lambda^2-2\lambda+t$& $(a,b,a,b)$\\ \hline
2C & $\lambda^2-2\lambda t +t$& $(b,a,a,b)$\\ \hline\hline
3A & $\lambda^4-6\lambda^2 t+4\lambda t +4\lambda t^2 -3t^2$& $(a,9a,a,a)$\\ \hline
3B & $3\lambda^4-4\lambda^3-4\lambda^3 t+6\lambda^2 t-t^2$& $(9a,a,a,a)$\\ \hline
3C & $\lambda^4-4\lambda^3+6t\lambda^2-4t^2\lambda+t^2$& $(a,a,9a,a)$\\ \hline
3D & $\lambda^4-4t\lambda^3+6t\lambda^2-4t\lambda+t^2$& $(a,a,a,9a)$\\ \hline\hline
4A & $\lambda^4-2t\lambda^3-2\lambda^3+6t\lambda^2$ &$(a,\frac{1}{8},a,a)$\\ \hline
   & $-2t^2\lambda-2t\lambda+t^3-t^2+t$& \\ \hline
4B & $\lambda^4-2t\lambda^3+2t^2\lambda-t^3$& $(a,a,\frac{1}{8},a)$\\ \hline
4C & $\lambda^4(t^2-t+1)-2\lambda^3t(t+1)+6t^2\lambda^2$& $(\frac{1}{8},a,a,a)$\\
   & $-2\lambda t^2(t+1)+t^3$& \\ \hline
4D & $\lambda^4-2\lambda^3+2t\lambda-t$& $(a,a,a,\frac{1}{18})$ \\ \hline\hline
5A & $-2\lambda^3+3t\lambda^2+3\lambda^2-6t\lambda+t^2+t$& $(4a,\frac{1}{8},a,a)$\\ \hline
5B & $\lambda^3-3\lambda^2+3t\lambda-2t^2+t$& $(a,\frac{1}{18},4a,a)$\\ \hline
5C & $\lambda^3-3t\lambda^2+3t\lambda+t^2-2t$& $(a,\frac{1}{18},a,4a)$\\ \hline
5D & $2\lambda^3-3t\lambda^2+t^2$ & $(4a,a,\frac{1}{18},a)$\\ \hline
5E & $\lambda^3-3t\lambda+2t^2$& $(a,4a,\frac{1}{18},a)$\\ \hline
5F & $\lambda^3-3t\lambda^2+3t\lambda-t^2$& $(a,a,\frac{1}{18},4a)$\\ \hline
5G & $\lambda^3 (2-t)-3t\lambda^2+3t^2\lambda-t^2$& $(\frac{1}{18},a,4a,a)$\\ \hline
5H & $\lambda^3 (t+1)-6t\lambda^2+3t(t+1)\lambda-2t^2$& $(\frac{1}{18},4a,a,a)$\\ \hline
5I & $(1-2t)\lambda^3+3t\lambda^2-3t\lambda+t^2$ & $(\frac{1}{18},a,a,4a)$\\ \hline
5J & $\lambda^3-3\lambda^2+3t\lambda-t$ & $(a,a,4a,\frac{1}{18})$\\ \hline
5K & $\lambda^3-3t\lambda+2t$& $(a,4a,a,\frac{1}{18})$\\ \hline
5L & $\lambda^3-3\lambda^2+t$& $(4a,a,a,\frac{1}{18})$\\ \hline\hline
\end{tabular}
\end{center}

\caption{List of all algebraic solutions satisfying families of $ \mathbf{PVI_\alpha}$  equations }
\label{table1}
\end{table}
%
%\begin{table}
%\input{table2.tex}
%
%\caption{List of all families of $\mathbf{PVI_\alpha}$  equations having a common algebraic solution}
%\label{table2}
%\end{table}
%\begin{figure}
%\input{fig1.pic}
%\caption{The families of $\mathbf{PVI_\alpha}$ equations whose solution is $2A$, $2B$, or $2C$}
%\label{fig1}
%\end{figure}
%\begin{figure}
%\input{fig2.pic}
%\caption{The families of $\mathbf{PVI_\alpha}$ equations whose solution is $3A$, $3B$,..., $4C$, or $4D$}
%\label{fig2}
%\end{figure}
%\begin{figure}
%\input{fig3.pic}
%\caption{The families of $\mathbf{PVI_\alpha}$ equations whose solution is $5A$, $5B$, ..., $5K$, or $5L$}
%\label{fig3}
%\end{figure}
{\bf Outline of the Proof.}  Denote by $\Gamma_\beta$ the  compactified and normalized algebraic curve defined
by (\ref{1}), with affine model
\begin{equation}\label{gb}
\Gamma_\beta^{aff} = \{(\lambda,t)\in \C^2: N(\lambda,t) =0 \}
\end{equation}
 where
\begin{eqnarray*}
  N(\lambda,t) &=& \beta_0 \lambda ^2(\lambda - 1)^2 (\lambda - t)^2 - \beta_1 t (\lambda - 1)^2 (\lambda - t)^2\\
    & +& \beta_2 (t-1)\lambda^2 (\lambda - t)^2
  - \beta_3 t(t - 1)\lambda ^2 (\lambda - 1)^2 = 0 .
\end{eqnarray*}
In the case when $\Gamma_\beta$ is irreducible, the relation $\{N(\lambda,t)=0\}$ defines an algebraic function
$\lambda(t)$. If this function were a solution of some $\mathbf{PVI_{\alpha}}$ equation, then the only
ramification points of $\lambda(t)$ would be at $t=0,1,\infty$ (because $\mathbf{PVI_{\alpha}}$ satisfies the so
called Painlev\'{e} property \cite{8}). Equivalently, the pair $(\Gamma_\beta,t)$ is a Belyi pair, which means that
the only possible critical values of the map
\begin{equation}
\label{projection} \pi: \Gamma_\beta \rightarrow \C\P^1: (\lambda,t) \rightarrow t
 \end{equation}
are $0$, $1$ or $\infty$. This means also that if $\Delta(t)$ is the discriminant of $N(\lambda,t)$ with respect
to $\lambda$, then it  is a polynomial whose only roots are at $t=0$ and $t=1$. A direct computation shows that
this is impossible. The more difficult case is when $N(\lambda,t)$ is reducible over $\C$. Then $\Gamma_\beta$
defines several algebraic functions and we have to apply the above to each of them. Finally we have to check
whether the obtained function is actually a solution of some $\mathbf{PVI_{\alpha}}$ equation. To check whether
a given polynomial $N(\lambda,t)$ is reducible over $\C$ is a difficult task in general. We shall make use of
the action of the symmetric group $\mathcal{S}_4$ (see section \ref{s4}) on the set of curves $\Gamma_\beta$,
parameterized by $\beta \in \C\P^3$.

It turns out that, first,  curves $\Gamma_\beta$ with a trivial stabilizer under the action of $\mathcal{S}_4$
can not produce a solution of $\mathbf{PVI_{\alpha}}$. The stabilizer of a curve acts on it as a group of
automorphisms (symmetries) which imposes additional restrictions on $\beta$.

The second ingredient of the proof is the study of the Puiseux expansion of $\lambda(t)$ in a neighborhood of
$t=0,1,\infty$ (section \ref{pui}). These expansions depend on the stabilizer of $\Gamma_\beta$ only and imply
the possible topological types of the solution $\lambda(t)$. Equivalently, to each solution $\lambda(t)$ we
associate a Belyi pair and the Puiseux expansions determine their possible \textit{dessin d'enfant}. The
algebraic functions which we obtain in this way are, \textit{a posteriori}, the solutions of the
$\mathbf{PVI_{\alpha}}$ presented in Table \ref{table1}.
\\
 \textbf{Proof of Theorem \ref{t7}.}
\subsection{The action of $\mathcal{S}_4$}
\label{s4}
 The set of automorphisms of the projective line $\C\P^1$ which send four distinct points
$(0,1,t,\infty)$ to the points $(0,1,\tilde{t},\infty)$ ($\tilde{t}=\tilde{t}(t)$ is uniquely defined) form a
group isomorphic to $\mathcal{S}_4$  generated by the transpositions
\begin{equation}\label{action1}
x^1:  s \mapsto 1-s, x^2:  s \mapsto \frac{1}{s}, x^3:  s \mapsto \frac{t-s}{t-1} .
\end{equation}
Each $x^i$ sends an isomonodromic family of Fuchsian systems with singular points at $0,1,t,\infty$ to an
isomonodromic family of such systems with singular points at $0,1,\tilde{t},\infty$. Therefore $x^i$ induce an
action of $\mathcal{S}_4$ on the set of $\mathbf{PVI_{\alpha}}$ equations, and hence on the set of curves
$\Gamma_\beta$. Explicitly we have
\begin{equation}\label{action2}
x^i: \Gamma_\beta \rightarrow \Gamma_{x^i_*(\beta)}: (\lambda,t) \rightarrow (x^i(\lambda),x^i(t) ), i=1,2
\end{equation}
and
\begin{equation}\label{action23}
x^3: \Gamma_\beta \rightarrow \Gamma_{x^3_*(\beta)}: (\lambda,t) \rightarrow (x^3(\lambda),x^3(0) )=
(\frac{t-\lambda}{t-1}, \frac{t}{t-1})
\end{equation}
where
\begin{equation}\label{action3}
\left\{
\begin{array}{lrc}
  x^1_* : (\beta_0,\beta_1,\beta_2,\beta_3) & \rightarrow &(\beta_0,\beta_2,\beta_1,\beta_3)\\
  x^2_* : (\beta_0,\beta_1,\beta_2,\beta_3) &\rightarrow &(\beta_1,\beta_0,\beta_2,\beta_3) \\
  x^3_* : (\beta_0,\beta_1,\beta_2,\beta_3) &\rightarrow &(\beta_0,\beta_3,\beta_2,\beta_1)
\end{array}\right.
\end{equation}
which is the standard representation of $\mathcal{S}_4$ on $\C^4$ (upon identifying $\infty,0,1,t$ to
$\beta_0,\beta_1,\beta_2,\beta_3$ respectively). The proof of the above facts is straightforward, see \cite{8}
for details.
\subsection{The
topological type of the projection $\Gamma_\beta \rightarrow \C\P^1$
 in a neighborhood of the pre-image
of $t=0,1,\infty$} \label{pui}
 Let $\Gamma_\beta$ be the compactified and normalized curve  defined by (\ref{gb}) (it is a
disjoint union of Riemann surfaces). In this section we  determine the topological type of the projection
(\ref{projection})
$$
\Gamma_\beta \rightarrow \C\P^1
$$
in a neighborhood of the pre-image of $t=0,1,\infty$ in $\Gamma_\beta$. In the projective space $\C\P^3$ with
coordinates $[\beta_0:\beta_1:\beta_2:\beta_3]$ consider the complex polyhedron $W$ formed by the ten planes
(2-faces)
\begin{equation}\label{relations}
W = \cup_{i\neq j}\{\beta_i= \beta_j\} \cup_k \{\beta_k=0\} .
\end{equation}
It has also $45$ 1-faces (projective lines) and $120$ 0-faces (points). We shall see in the process of the proof
that the topological type of the projection in a neighborhood of the pre-image of $t=0,1,\infty$ is one and the
same when $\beta$ belongs to a given $i$-face, but does not belong to any other $j$-face with $j<i$. For this
reason we shall use, until the end of this paper, the following convention. \emph{When we say that a point
$\beta$ belongs to a given face (satisfies some set of relations (\ref{relations})), then this will mean that it
does not belong to any other face of smaller dimension (does not satisfy any other relation from the list
(\ref{relations}) ). }
 The topological type in a neighborhood of the pre-image of
any point is determined by a partition of the degree of the map which is $6$. Thus a partition $(1+1+1+1+2)$
means that we have $5$ pre-images, and that the multiplicity of $\pi$ at each pre-image  is $1,1,1,1,2$
respectively. Similarly, a partition $(1+1+2+2)$ means that have $4$ pre-images  with multiplicities $1,1,2,2$
respectively etc. To formulate the result we note that the symmetric group $\mathcal{S}_4$ acts on the
polyhedron $W$ by its standard representation (\ref{action3}), as well on the set of curves $\Gamma_\beta$ by
(\ref{action2}). The subgroup $\mathcal{S}_3$ generated by $x^1, x^2$ permutes the ramification points
$0,1,\infty$ according to (\ref{action1}) without changing the topological type of the projection $\pi$ over
each of these points.
\begin{Pro}
The topological type of the projection (\ref{projection}) in a neighborhood of the pre-image of $t=0,1,\infty$
is one and the same when $\beta$ belongs to a given face of the polyhedron $W$ or it does not belong to $W$.
This topological type is shown on Table \ref{table3} (one representative for each orbit of $\mathcal{S}_3=
<x^1,x^2>$)
\end{Pro}
\begin{proof}[Proof] The bi-rational transformations $x^1, x^2$ defined by (\ref{action2}) are compatible with the
projection $\pi$ and permute the points $t=0,1,\infty$.
 Therefore it suffices to consider the pre-image of $0$.
\begin{figure}
\input{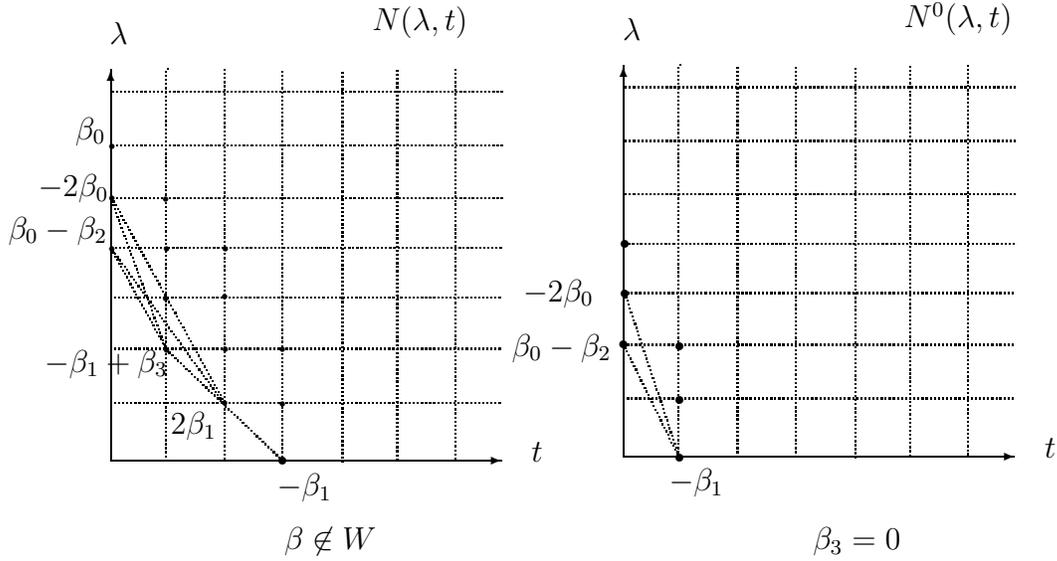}
  \caption{The Newton polygon of $N(\lambda,t)$ and $N^0(\lambda,t)$}
\label{fig4}
\end{figure}
Let us consider in detail the "generic" case, when $\beta\not\in W$. It follows from the Newton polygon of
$N(\lambda,t)$,  shown on fig. \ref{fig4}, that there are at least three Puiseux series in a neighborhood of
$(0,0)$ (for the terminology see for instance \cite[Kirwan]{kirwan}). The first two correspond to the line
segment $[(3,0),(1,2)]$ and have non-equivalent leading terms
$$
\lambda = c_1 t + \dots, \lambda = c_2 t + \dots
$$
where
$$
(\beta_3-\beta_1) c^2_{1,2} + 2 \beta_1 c_{1,2} - \beta_1 = 0
$$
provided that
$$
 \beta_1\neq \beta_3, \beta_1^2 + \beta_1 (\beta_3-\beta_1) \neq 0 .
$$
The third one corresponds to the line segment $[(1,2),(0,4)]$ and has leading term
$$
\lambda = c_3 t^{1/2} + \dots
$$
where
$$
(\beta_0-\beta_2) c_3^2 + \beta_3 - \beta_1 =0,
$$
provided that
$$
\beta_0 \neq \beta_2, \beta_1\neq \beta_3 .
$$
 Taking into consideration that
$$
N(\lambda,0)= \lambda^4 ( \beta_0 \lambda^2 -2 \beta_0 \lambda + \beta_0 - \beta_2)
$$
we conclude that we have at least five pre-images of multiplicities at least $1,1,2,1,1$ respectively. As the
degree of the map $\pi$ is six, then its topological type is exactly
 $(1+1+1+1+2)$. The topological type of the projection $\pi$ over $0$ and $1$ is obtained by acting with
the group $\mathcal{S}_3$ generated by $x^1$, $x^2$.

In a similar way one verifies that when $\beta_0 = \beta_2$, or $\beta_1= \beta_3$ the multiplicities are
$(1+1+1+3)$. If
 $ \beta_0 = \beta_2$ and $ \beta_1= \beta_3 $ the multiplicities are $(1+1+2+2)$. The case $ \beta_0 =  \beta_1= \beta_3 $
 is the same as $ \beta_0 = \beta_2$ and the multiplicity is $(1+1+1+3)$. The case
 $ \beta_0 = \beta_1 = \beta_2= \beta_3 $ is of the same type as $ \beta_0 = \beta_2$ and $ \beta_1= \beta_3 $.
The multiplicities of $\pi$ over $1$ and $\infty$ are obtained as before. This completes the study of faces of
$W$ for which $\beta_i\neq 0$.
In  the case $\beta_3=0$ we consider the curve
$$ \Gamma_\beta^0 = \{(\lambda,t)\in \C^2: \beta_0 - \beta_1 \frac{t}{\lambda ^2} + \beta_2
\frac{t - 1}{(\lambda - 1)^2}
   = 0, \lambda\neq 0,1
   \} .
  $$
The polynomial $N(\lambda,t)$ is replaced by
$$N^0(\lambda,t)= \beta_0 \lambda ^2(\lambda - 1)^2  - \beta_1 t (\lambda - 1)^2
    + \beta_2 (t-1)\lambda^2
    $$
whose  Newton polygon
  is shown    on fig. \ref{fig4}.
It follows that there is at least one Puiseux expansion with leading term
$$
\lambda = c t^{1/2} + \dots,  \beta_1 + (\beta_2- \beta_0) c^2=0
$$
provided that $\beta_1 \neq 0$, $ \beta_2\neq \beta_0$. As
$$
N^0(\lambda,0)= \lambda^2 (\beta_0 (\lambda-1)^2 - \beta_2)
$$
then $t=0$ has at least three pre-images, provided that $\beta_0 \beta_2 \neq 0$. We
conclude that $t=0$ has exactly three pre-images with multiplicities $2,1,1$ respectively, provided that $\beta$
belongs to the 2-face $\beta_3=0$. The remaining 1-faces and 0-faces are studied in the same way. The result is
summarized
 on Table \ref{table3}.
 It worth noting that in all cases the computing  of the
leading term of the Puiseux expansion suffices to deduce the result.
\end{proof}

We conclude this section by the
following elementary claim which will be often useful in the computations
\begin{table}
  \centering
  \begin{tabular}{|c|c|c|c|c|}
    \hline
     {\rm face of $W$}&stabilizer & $t=0$ & $t=1$ & $t=\infty$ \\
    \hline
    $\beta_i\neq \beta_j $ &$S_4$ & (1+1+1+1+2) & (1+1+1+1+2) & (1+1+1+1+2) \\
    $\beta_0=\beta_2$ & $S_2 \times S_2$ &(1+1+1+3) & (1+1+1+1+2) & (1+1+1+1+2) \\
    $\beta_0=\beta_2, \beta_1=\beta_3$& $D_4$ &(1+1+2+2) & (1+1+1+1+2) & (1+1+1+1+2) \\
    $\beta_0=\beta_1=\beta_2$ & $S_3$& (1+1+1+3) & (1+1+1+3) & (1+1+1+3) \\
    $\beta_0=\beta_1=\beta_2=\beta_3 $&$ S_4$ & (1+1+2+2) & (1+1+2+2) & (1+1+2+2) \\
    $ \beta_3=0 $&$S_3$ & (1+1+2) & (1+1+2) & (1+1+2) \\
    $\beta_0=\beta_2, \beta_3=0 $& $S_2$&(1+3) & (1+1+2) & (1+1+2) \\
    $\beta_0=\beta_1=\beta_2, \beta_3=0 $& $S_3$ &(1+3) & (1+3) & (1+3) \\
    $\beta_2= \beta_3=0$               & $S_2\times S_2$  & (2) & (1+1) & (2)\\
    $\beta_2= \beta_3=0, \beta_0=\beta_1$               & $S_2\times S_2$  & (2) & (1+1) & (2)\\
    \hline
  \end{tabular}
 \caption{Multiplicity of $\pi$ at the pre-images of $t=0,1,\infty$. }
 \label{table3}
\end{table}
\begin{Pro}
\label{n1}
 Let $N_1(\lambda,t)$ be a polynomial of non-zero degree with respect to $\lambda$ and of non-zero
degree with respect to $t$, which divides $N(\lambda,t)$, and $\beta_1 \beta_2 \beta_3\neq 0$.
 Then
$$
N_1(0,t)= c_0 t^{n_0}, c_0\neq 0, 1\leq n_0 \leq 3, N_1(1,t)= c_1 (t-1)^{n_1}, c_1\neq 0, 1\leq n_1 \leq 3
$$
and
$$
N_1(t,t)= c_2 t^{m_0} (t-1)^{m_1}, c_2\neq 0, 1\leq m_0, 1\leq m_1, m_0+m_1 \leq 3 .
$$
\end{Pro}
{\bf Proof.} We have $N(0,t)= -\beta_1 t^3$. For a fixed $\lambda=c\sim 0$ the polynomial $N(c,t)\in \C[t]$ has
exactly three roots which tend to zero when $c$ tends to zero. Therefore the polynomial $N_1(c,t)\in \C[t]$ has
at least one and at most three roots which tends to zero when $c$ tends to zero, which proves the claim
concerning $N_1(0,t)$. The claim concerning $N_1(1,t)$ is proved in the same way. As $N_1(0,0)=N_1(1,1)=0$ then
$N_1(t,t)$ is divided by $t(t-1)$ but also divides $N(t,t)= -\beta_3 t^3(t-1)^3$. $\Box$

We are ready to compute the solutions of $ \mathbf{PVI_\alpha}$ corresponding to the faces of $W$. Let $\Gamma$
be the Riemann surface of an irreducible  component of $\Gamma_\beta$, which defines a solution of some $
\mathbf{PVI_\alpha}$ equation. Then the only ramification points of the induced map
\begin{equation}\label{bpi}
\pi :\Gamma \rightarrow \mathbb{C}\mathbb{P}^1: (\lambda,t) \rightarrow t
\end{equation}
are at $0,1,\infty$, and $\Gamma$ is connected. The pair $(\Gamma,\pi)$ is called a Belyi pair, to which we
associate a \emph{dessins d'enfant}, which is the graph obtained as a pre-image of the segment $[0,1]$ under the
map $\pi$. The degree of the dessin is the degree of $\pi$ (see \cite{sch}). The dessin d'enfant will be useful
when describing the topological type of the projection $\pi$.

\subsection{The case $\beta\not\in W$}
\label{notb}
\begin{figure}
  \input{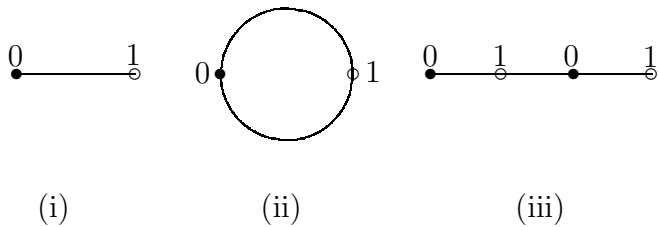}
  \caption{Dessins of degree $1,2$ and $3$}
\label{fig5}
\end{figure}
 We suppose that $\lambda(t)$ is an algebraic function, such that $N(\lambda(t),t)\equiv 0$, and
consider the corresponding Belyi pair $(\Gamma,\pi)$, (\ref{bpi}). Let $\{ \lambda_1,\dots,\lambda_d\}=
\pi^{-1}(t_0)$ where $t_0\neq 0,1,\infty$. The loops originating from $t_0$ and going clockwise once around
$0,1$ and $\infty$ induce permutations $\sigma_0,\sigma_1$ and $\sigma_\infty$ of the points
$\lambda_1,\dots,\lambda_d$, such that $\sigma_0\sigma_1\sigma_\infty= 1$. According to Table \ref{table3}
$\sigma_0,\sigma_1$ and $\sigma_\infty$ are transpositions, unless one of them is the identity permutation. In
the former case $ (\sigma_0 \sigma_1)^2=1$ and hence $\sigma_0\sigma_1 = \sigma_1\sigma_0$. Thus $\sigma_\infty$
is a product of two disjoint transpositions, which contradicts to Table \ref{table3}. On the other, if one of
the permutations $\sigma_0,\sigma_1,\sigma_\infty$ is the identity,  the group generated by them
 is either $\mathbf{Z}_2$ or is trivial. This shows that the degree of $\pi$ is
either two (because the covering (\ref{bpi}) is connected), or one. The corresponding dessin d'enfants are
 shown on fig.\ref{fig5} (i) and (ii) (in particular
$N(\lambda,t) \in \C[\lambda,t]$ is reducible). If the dessin is of degree one, then the solution is defined as
$\lambda=P(t)$, where $P$ is a polynomial (the coefficient of $\lambda^6$ in the polynomial $N(\lambda,t)$ is
$\beta_0\neq 0$). By Proposition \ref{n1} we conclude that either $\lambda=t$, or $\lambda=t^2$ or
$\lambda=t^3$. But $\lambda-t$ , $\lambda-t^2$, $\lambda-t^3$ can not divide $N(\lambda,t)$, provided that
$\beta_i\neq 0$. If the dessin is of degree two, then $\lambda(t)$ is defined by $\lambda^2+2 p(t)\lambda +
q(t)=0$. The functions $p,q$ are polynomials in $t$, because the coefficient of $\lambda^6$ in the polynomial
$N(\lambda,t)$ is $\beta_0\neq 0$. Further, we may suppose (acting with an appropriate symmetry $x^i$ on
$\Gamma$, see section \ref{s4}) that
 $\lambda(t)$ is ramified over $0$ and $\infty$ only.
 By Proposition \ref{n1} $q(t)$ is a non-constant polynomial which divides $t^3$. As $p(t)^2-q(t)$ is a non-constant
 monomial of odd degree, then $p(t)\equiv 0$. Proposition \ref{n1} implies that we have either $\lambda^2=t$ or
$\lambda^2=t^3$. The polynomial $\lambda^2-t^3$ can not divide $N(\lambda,t)$ while $\lambda^2-t$  divides
$N(\lambda,t)$ if and only if $\beta_0=\beta_1$ and $\beta_2=\beta_3$ (this case is excluded, as $\beta\not\in
W$). \textit{The curve $\Gamma_\beta$ does not define  a solution.}
\subsection{The face $\beta_1 = \beta_2$}
\label{b1b2}
 The possible  dessins d'enfant are determined as above. Namely, when one of the permutations
$\sigma_0,\sigma_1,\sigma_\infty$ is identity, the dessin is of degree one or two. Up to a symmetry it is
equivalent to the one shown on fig. \ref{fig5} (i) or (ii). Reasoning as in the case $\beta \not \in W$ we
conclude that $\Gamma_\beta$ does not define a solution.

If, on the other hand, $\sigma_0$, $\sigma_1$ are non-trivial transpositions we have one more case compared to
section \ref{notb}: $\sigma_3$ is cyclic of order three, and $\sigma_0$, $\sigma_1$ are non-disjoined
permutations. Taking into account that the covering (\ref{bpi}) is connected, we conclude that the dessin is of
degree three. Up to a symmetry, it is shown on \ref{fig5} (iii) and $\lambda(t)$ satisfies
$N_1(\lambda(t),t)\equiv 0$ where $N_1$ is an irreducible polynomial of degree three in $\lambda$ dividing the
polynomial $N(\lambda,t)$, defined after formula (\ref{gb}). We denote
 $$
 N=N_1 N_2, \Gamma_1^{aff}=\{N_1(\lambda,t)=0\}, \Gamma_2^{aff}=\{N_2(\lambda,t)=0\},\Gamma_\beta^{aff}= \Gamma_1^{aff}\cup
 \Gamma_2^{aff} .
 $$
 As before, let $\Gamma_1, \Gamma_2, \Gamma_\beta$ be the corresponding compactified and normalized curves.
 The symmetry $x^1$ is an automorphism of $\Gamma_\beta$ and hence it is either an automorphism of
the curves $\Gamma_1$ and $\Gamma_2$, or it permutes these curves (we used that $\Gamma_1$ is irreducible).
Suppose first that $x^1$ is an automorphism of $\Gamma_1$. Then the rational function
$$
 \frac{N_1(\lambda,t)}{\lambda(\lambda-1)(\lambda-t)}= 1+\frac{A_1}{\lambda}+
\frac{B_1}{\lambda-1}+\frac{C_1}{\lambda-t}
$$
is invariant under the action of $x^1$ too. Here $A_1,B_1,C_1$ are polynomials in $t$ which divide $t, t-1$ and
$t(t-\lambda)$ respectively (see (\ref{1})), and hence we have
$$
A_1(1-t)=-B_1(t), B_1(1-t) = -A_1(t), C_1(1-t) = -C_1(t) .
$$
Similarly, if
$$
\frac{N_2(\lambda,t)}{\lambda(\lambda-1)(\lambda-t)}= 1+\frac{A_2}{\lambda}+
\frac{B_2}{\lambda-1}+\frac{C_2}{\lambda-t}
$$
then
$$
A_2(1-t)=-B_2(t), B_2(1-t) = -A_2(t), C_2(1-t) = -C_2(t) .
$$
We conclude that $C_1(t)= c_1 t(t-1)$, $C_2(t)= c_2 t(t-1)$, which contradicts to $C_1C_2= -
\beta_3t(t-1)/\beta_1$, $\beta_1,\beta_3\neq 0$.
Suppose now that the map $x^1$ exchanges the curves $\Gamma_1$ and $\Gamma_2$. Then we have
$$
A_1(1-t)=-B_2(t), B_1(1-t) = -B_2(t), C_1(1-t) = -C_2(t) .
$$
The polynomial $N(\lambda,t)$ is of degree three with respect to $t$ and
$$
A_1A_2 = - \frac{\beta_1}{\beta_0}t, B_1B_2 = -\frac{\beta_2}{ \beta_0}(t-1), C_1C_2 = -
\frac{\beta_3}{\beta_0}t(t-1) .
$$
Therefore without loss of generality we may suppose that $C_1(t) = c_1 t, C_2= c_1 (t-1)$ and $A_1(t)= a_1 t,
B_2=b_2 (t-1)$ or $A_2(t)= a_2 t, B_1=b_1 (t-1)$, where $a_i,b_j\neq 0$. In both cases the polynomials
$N_1(\lambda,t)$, $N_2(\lambda,t)$ are of degree two in $t$, in contradiction to the fact that the degree of
$N(\lambda,t)$ with respect to $t$ is three. \textit{We conclude that the curve $\Gamma_\beta$ does not define a
solution.}
\subsection{The face $\beta_0 = \beta_2$, $\beta_1 = \beta_3$.}
We have the identity
$$
N(\lambda,t) = (\lambda^2- 2 \lambda +t)(\beta_0\lambda^2(\lambda-t)^2 - \beta_1 t^2 (\lambda -1)^2) .
$$
Indeed,
\begin{equation}\label{1B} \lambda^2- 2 \lambda +t=0
\end{equation}
 defines a solution of $\mathbf{PVI_{\alpha}}$, e.g. \cite{doran}, Table 2,
solution $2B$. Its dessin is equivalent to the one on fig. \ref{fig5} (ii). The function $\lambda(t)$ defined by
$$\lambda(\lambda-t) - c t (\lambda
-1)=0, c= \pm \sqrt{\frac{\beta_1}{\beta_0}}
$$
is ramified over $0,1,\infty$ only provided that $c=0,\pm 1$. This is, however, impossible as $\beta_0\neq
\beta_1, \beta_i\neq 0$.
\textit{The curve $\Gamma_\beta$  defines  the solution (\ref{1B}).}
\subsection{The face $\beta_0=\beta_1 = \beta_2$.}
According to Table \ref{table3}  each of the permutations $\sigma_0,\sigma_1,\sigma_\infty$ is either the
identity, or is a cycle of length three.

 If one of the permutations $\sigma_0,\sigma_1,\sigma_\infty$ is the identity then the group
generated by $\sigma_0,\sigma_1,\sigma_\infty$ is either $\mathbf{Z}_3$ or the trivial one $\{\mathbf{1}\}$, and
hence the degree of the corresponding dessin is  one or three. The case of degree one does not lead to a
solution (see section \ref{notb}). The case of degree three is studied as in section \ref{b1b2} and does not
lead to a solution too (provided that $\beta_i \neq 0$).

If neither of the permutations $\sigma_0,\sigma_1,\sigma_\infty$ is the identity, then they are disjoint
three-cycles. As the symmetric group $\mathcal{S}_3$ contains only two three-cycles we conclude that the degree
of the projection $\pi$ is at least four. Suppose that $\lambda(t)$ is defined by  the  polynomial $N_1$, $
N_1(\lambda(t),t)\equiv 0$,
 where $N_1\in \C[\lambda,t]$ is
irreducible of degree four in $\lambda$. Then $x^1, x^2$ are automorphisms of $$\Gamma_1 = \{ N_1(\lambda,t)=0
\}.$$ It follows that the curve
$$\Gamma_2 = \{ N_2(\lambda,t)=0 \}$$
 defined by the polynomial $N_2=N/N_1$ is also invariant. We have
$$
\frac{N_2(\lambda,t)}{\lambda(\lambda-1)}= 1+ \frac{A}{\lambda} + \frac{B}{\lambda-1}
$$
where $A,B$ are polynomials in $t$ of degree at most three. The $x^{1,2}$ invariance of the above expression
implies
$$
A(t)A(1/t)=1, A(1/t) B(t) = - B(1/t), B(t)=-A(1-t) .
$$
with solutions
$$
A(t)= t, B(t)= t-1; A(t)= t^3, B(t)= (t-1)^3; A(t)= -t^2, B(t)= (t-1)^2 .
$$
The case $A(t)= t^3, B(t)= (t-1)^3$ does not lead to a solution as $N_1(\lambda,t)$ does  depend on $t$. The
case $A(t)= -t^2, B(t)= (t-1)^2$ implies $N_2(\lambda,t)= (\lambda-t)^2$, and hence $\beta_3=0$. Finally, in the
case $A(t)= t, B(t)= t-1$ we have $N_2(\lambda,t)= \lambda^2-2\lambda + 2\lambda t - t$ (which does not define a
solution). The condition that $N_2(\lambda,t)$ divides $N(\lambda,t)$ leads to $\beta_3= 9 \beta_0$ and we get
$$
N(\lambda,t)= \left( {{\it \lambda}}^{2}-2\,{\it \lambda}+2\,{\it \lambda}\,t-t \right)  \left( {t}^{2}-4\,{{\it
\lambda}}^{3}t+6\,{{\it \lambda}}^{2}t-4\,{\it \lambda}\,t+{{\it \lambda}}^ {4} \right) .
$$
The function $\lambda(t)$ defined by
\begin{equation}\label{3D}
{t}^{2}-4\,{{\it \lambda}}^{3}t+6\,{{\it \lambda}}^{2}t-4\,{\it \lambda}\,t+{{\it \lambda}}^ {4}=0
\end{equation}
 is indeed a solution of $\mathbf{PVI_{\alpha}}$, e.g. \cite{doran}, Table 2,
solution $3D$. In the case when the dessin corresponding to $\lambda(t)$ is of degree five we conclude that the
polynomial $N_2(\lambda,t)$ is linear in $\lambda$. By Proposition \ref{n1} we get $N_2(\lambda,t)=\lambda-t^2$
which implies $\beta_1=0$. To resume, \textit{the curve $\Gamma_\beta$  defines  a solution, provided that
$\beta_0=\beta_1 = \beta_2=\beta_3/9$.}
\subsection{The face $\beta_0=\beta_1 = \beta_2=\beta_3$.}
We have
$$
N(\lambda,t)=\beta_0 (\lambda^2- 2 \lambda +t)(\lambda^2- 2 \lambda t +t)(\lambda^2 -t)
$$
and the three algebraic functions defined by $N(\lambda,t)=0$ are solutions of suitable $\mathbf{PVI_{\alpha}}$
equations, e.g. \cite{doran}, Table 2, solutions $2B,2C,2A$ respectively.
\textit{The curve $\Gamma_\beta$  defines  three solutions}
\subsection{The face $\beta_3 =0$.}
Recall that in this case $N(\lambda,t)= (\lambda-t)^2N^0(\lambda,t)$ where
$$
N_0(\lambda,t)= \beta_0 \lambda ^2(\lambda - 1)^2  - \beta_1 t (\lambda - 1)^2
    + \beta_2 (t-1)\lambda^2 .
    $$
 The same arguments as in section \ref{notb} show that the corresponding dessin d'enfant
is of degree one or two.

If the degree is one, then the solution is $\lambda=P(t)$ for some non-constant polynomial $P$. Therefore
$N^0(P(t),t) \not\equiv 0$ and $P(t)$ can not be a solution.

If the degree is two, then $\lambda(t)$ has exactly two ramification points. Without loss of generality we
suppose that these points a $0$ and $\infty$, and as in section \ref{notb} we conclude that $\lambda(t)$ is
defined by $\lambda^2+2 p(t)\lambda + q(t)=0$ for some  $p,q\in \mathbb{C}[t]$. The polynomial $\lambda^2+2
p(t)\lambda + q(t)$ divides
$$
N_0(\lambda,t)= t(\beta_2 \lambda^2 - \beta_1  (\lambda - 1)^2) +  \beta_0 \lambda ^2(\lambda - 1)^2 -\beta_2
\lambda^2
$$
and hence $p(t)=c_1$ and $q(t)=c_2 t$ for some constants $c_1, c_2$. Without loss of generality we suppose that
the ramification points of $\lambda(t)$ are   $0$ and $\infty$ and hence the discriminant $4(p^2-q)$ is a power
of $t$.  This implies that $c_1=0$. Finally, a direct computation shows that the identity
$N^0(\sqrt{-c_2t},t)\equiv 0$ implies $\beta_2=0$ which is not true. \textit{The curve $\Gamma_\beta$ does not
define  a solution.}
\subsection{The face $\beta_3 =0,\beta_0=\beta_2$.}
It is easier to analyze the face $\beta_3 =0,\beta_1=\beta_2$, which is equivalent to $\beta_3
=0,\beta_0=\beta_2$ after applying the transformation $x^2$. Suppose for a moment that  $\beta_3
=0,\beta_1=\beta_2$. The dessin is of degree at most three, and hence $N^0(\lambda,t)$ is reducible.
 It follows that
 $\lambda-c$ divides $N^0(\lambda,t)$ for some constant $c$.
 As $N^0(\lambda,t)$ is linear in $t$,
then $\lambda-c$ is deduced from the coefficient of $t$ which equals to $\beta_1(1-2\lambda)$. Thus $c=1/2$ and
the condition that $1-2\lambda$ divides $N^0(\lambda,t)$ leads to $\beta_0=4\beta_1= 4\beta_2$, in which case
$$N^0(\lambda,t) = \beta_0
\left( 2\,{\it \lambda}-1 \right)  \left( 2\,{{\it \lambda}}^{3}-3\,{{\it \lambda}}^{2}+t \right) .
$$
The function $\lambda(t)$ defined by
 $\left( 2\,{{\it \lambda}}^{3}-3\,{{\it \lambda}}^{2}+t \right)=0$ is indeed a solution, see \cite{doran},
Table 2, solutions $5L$.  Applying the transformation $(x^2)^{-1}=x^2$ of section \ref{s4} we get the solution
(see \cite{doran}, Table 2, solution 5K
$$
 {{\it \lambda}}^{3}-3 \,{{\it \lambda}} t+ 2 \,t = 0
$$
defined by $\Gamma_\beta$ with $\beta_3 =0,\beta_1=4 \beta_0=4 \beta_2$.
\textit{The curve $\Gamma_\beta$  defines  a solution provided that $\beta_3 =0,\beta_1=4 \beta_0=4 \beta_2$.}
\subsection{The face $\beta_3 =0,\beta_0=\beta_1=\beta_2$.}
\textit{The polynomial $N^0(\lambda,t)$ is irreducible and defines a solution, see \cite{doran}, Table 2,
solution 4D. }
\subsection{The face $\beta_3 =0,\beta_2=0$.}
\emph{A direct computation shows that the relation $\beta_0 \lambda^2 - \beta_1 t=0$ defines a solution.}
\subsection{The face $\beta_3 =0,\beta_2=0,\beta_0=\beta_1$.}
\textit{The solution is $\lambda^2=t$.}

The results are summarized in Table \ref{table1}. Theorem \ref{t7} is proved.
\section{Algebraic solutions of  $\mathbf{PVI_\alpha}$ and Picard-Fuchs equations}
\label{section4} It was noted in the Remark after Theorem \ref{t7} that the families $1A, 1B,\dots ,1F$ are
Okamoto equivalent to the families $2A,2B,2C$. To this end we consider the remaining 23 families $2A-5L$, see
Table \ref{table1}. To each of them corresponds an affine plane or line in the parameter space $\mathbb{C}^4\{\alpha\}$
which, as we shall prove bellow, is generated by special points $\alpha$ of geometric origin, see Tables
\ref{table4} and \ref{table5}. Indeed, observe that exactly the same $23$ solutions $2A-5L$ were already
obtained by Doran, see Theorem \ref{t2} bellow,  by making use of deformations of elliptic surfaces with four
singular fibers. The corresponding special values of the parameter $\alpha$ are given in Table \ref{table4}. The
main result of this section is that exactly the same list of solutions can be obtained from deformations of
ramified covers of $\P^1$ with four ramification points. The corresponding values of the parameters $\alpha$ are
different and are shown on Table \ref{table5}, see Theorem \ref{t3}.

Recall that an elliptic surface  is a complex compact analytic surface $S$ with a projection $S \rightarrow
\P^1$, such that the general fiber $f^{-1}(z)=\Gamma_z$ is an elliptic curve. Two elliptic surfaces are
equivalent, if there is a bi-analytic map compatible with the projections, see \cite[Kodaira]{kodaira}.

We may suppose that the fiber $\Gamma_z$ is written in the Weierstrass form
$$
\Gamma_z = \{(x,y)\in \C^2: y^2= 4 x^3 - g_2(z) x - g_3(z) \}
$$
and consider the complete elliptic integrals of first and second kind
$$
\eta_1= \int_{\gamma(z)} \frac{dx}{y}, \ \eta_2= \int_{\gamma(z)} \frac{xdx}{y}
$$
where $\gamma(z) \subset \Gamma_z$ is a continuous family of closed loops (representing a locally constant
section $z\mapsto H_1(\Gamma_z,\Z)$ of the associated homology bundle). Then $\eta_1, \eta_2$ satisfy the
following Picard-Fuchs system (this goes back at least to \cite[Griffiths]{gri}, see \cite[Sasai]{sas})
\begin{equation}
\label{pf1}
\triangle(z) \frac{d}{dz}
\left(%
\begin{array}{c}
  \eta_1 \\
  \eta_2 \\
\end{array}%
\right) =
\left(%
\begin{array}{cc}
  -\frac{\triangle'_z}{12} & -\frac{3 \delta}{2} \\
  - \frac{g_2 \delta}{8} & \frac{\triangle'_z}{12}\\
\end{array}%
\right)
\left(%
\begin{array}{c}
  \eta_1 \\
  \eta_2 \\
\end{array}%
\right)
\end{equation}
where
$$
\triangle(g_2,g_3) = g_2^3- 27 g_3^2 .
$$
and
$$
\delta(z) = 3 g_3 \frac{d g_2}{dz} - 2 g_2 \frac{d g_3}{dz} .
$$
The singular points of the system correspond to singular fibers of the surface. The elliptic surfaces with four
singular fibers were classified by \cite[Herfurtner]{herfurtner} who obtained 50 distinct case, but only 5 of
them contain an additional parameter, see Table \ref{table6}. They lead to non-trivial isomonodromic
deformations of the above Picard-Fuchs system with four regular singular points. If we renormalize the singular
points to be $0,1,\infty,t$ then the zero $\lambda$ of $\delta(z)$, considered as a function in $t$ is a
solution of an appropriate $\mathbf{PVI_\alpha}$ equation, see \cite{okamoto} for details. The result is
summarized as follows
\begin{Thm}\cite[Theorem 3.13]{doran}
\label{t2}
 All algebraic solutions $(\lambda(t),\alpha)$
 of $\mathbf{PVI_\alpha}$ equation coming from moduli of elliptic surfaces with four singular fibers are shown
 on Table \ref{table1}. The corresponding values of $\alpha$ together with the stabilizer of
  the solution and the
$ \mathbf{PVI_\alpha}$ equation under the action of the symmetric group $S_4$
 are listed on Table \ref{table4}.
\end{Thm}
The Picard-Fuchs system (\ref{pf1}) has infinite monodromy group.
 We shall deduce similar families of Picard-Fuchs systems having a finite monodromy.
  Consider a ramified covering $\Gamma \rightarrow \P^1$ of degree three with
branching locus consisting of four points, where $\Gamma$ is a Riemann surface. We choose an affine model
$\Gamma^{aff}= \{(x,z)\in \mathbb{C}^2: f(x,z)=0 \} $ of $\Gamma$, where $f(x,z)= 4 x^3 - g_2 x - g_3$, and $ g_2=g_2(z),
g_3=g_3(z)$ are suitable polynomials. Moreover, without loss of generality, we suppose that the covering $\Gamma
\rightarrow \P^1$ is induced from the projection
\begin{equation}\label{cover}
    \{(x,z)\in \mathbb{C}^2: f(x,z)=0 \} \rightarrow \mathbb{C}: (x,z) \mapsto z .
\end{equation}
Let $x_1(z), x_2(z)$ be two distinct roots of $f$. Then $\gamma(z)=x_1(z) - x_2(z)$ is a 0-cycle of the fiber $\{x:
f(x,z)=0\}$ and the Abelian integrals above are replaced by
 the algebraic functions
$$\eta_1(z) = \int_{\gamma(z)}x = x_1(z)-x_2(z),  \eta_2(z)= \int_{\gamma(z)} x^2 = x_1^2(z)-x_2^2(z) .$$
A straightforward computation shows that $\eta_1, \eta_2$ satisfy the following Picard-Fuchs system (see
\cite{gm})
\begin{equation}
\label{pf2}
 \triangle(z) \frac{d}{dz}
\left(%
\begin{array}{c}
  \eta_1\\
  \eta_2 \\
\end{array}%
\right)
=
\left(%
\begin{array}{cc}
  \frac{\triangle'_z}{6} & -3 \delta \\
  - \frac{g_2 \delta}{2} & \frac{\triangle'_z}{3}\\
\end{array}%
\right)
\left(%
\begin{array}{c}
  \eta_1 \\
\eta_2\\
\end{array}%
\right) .
\end{equation}
As before, if we renormalize  the system (\ref{pf2}) to have singular points at $0,1,t, \infty$, then the root
$\lambda(t)$ of $ \delta(z)$ is a solution of a suitable  $ \mathbf{PVI_\alpha}$ equation, provided that
the deformation is isomonodromic. The last property holds, strictly speaking, in the case when the system is
non-resonant. In our case it holds too, because the deformation is  isoprincipal in the sense of
\cite{kats}. This can be also checked by a direct computation. Thus, if we consider the fibration (\ref{cover})
 and take for $g_2$, $g_3$ the expressions found by Herfutner, see Table
\ref{table6},  we get the same 23 algebraic solutions shown on Table \ref{table1}. The corresponding values for
$\alpha$ are of course different and are computed in Table \ref{table5}. In this way we prove
\begin{Thm}
\label{t3}
 The algebraic solutions $(\lambda(t),\alpha)$
 of $\mathbf{PVI_\alpha}$ equation coming from deformations of the covering (\ref{cover}) with
$g_2$, $g_3$ as on the Herfutner list, Table \ref{table6}, are shown
 on Table \ref{table1}. The corresponding values of $\alpha$ together with the stabilizer of
  the solution and the
$ \mathbf{PVI_\alpha}$ equation under the action of the symmetric group $S_4$
 are listed on Table \ref{table5}.
\end{Thm}
\begin{Rq}
As we already noted, the monodromy group of the Picard-Fuchs system (\ref{pf2}) is finite. Particular cases of
this system, in a more or less explicit way, were considered by many authors, e.g.  \cite[Boalch]{boalch},
\cite[Dubrovin-Mazzocco]{dub00}, \cite[Hitchin]{7,h7}, \cite[Kitaev]{kitaev}.
\end{Rq}

To this end we present, for convenience of the reader, one example of such computation. Using the Picard-Fuchs
system $(\ref{pf2})$, the Abelian integral of first kind $\eta_1$ satisfies the following equation
\begin{equation}
\label{pfeq}
p_0(z,a)\eta_1''+p_1(z,a)\eta_1'+p_2(z,a)\eta_1=0
\end{equation}
where
\begin{eqnarray*}
p_0(z,a)&=&144\delta\Delta^2\\
p_1(z,a)&=&144\Delta (\delta\frac{d\Delta}{dz}-\Delta\frac{d\delta}{dz})\\
p_2(z,a)&=&12\delta\frac{d^2\Delta}{dz^2}-216\delta^3g_2-12\Delta\frac{d\delta}{dz}\frac{d\Delta}{dz}-\delta(\frac{d\Delta}
{dz})^2 .
\end{eqnarray*}

Consider, for instance, the deformation $2$ from the  Herfurtner list (Table~\ref{table6})
\begin{table}

\begin{center}
\begin{tabular}{||c|c||}
\hline\hline
name & deformation \\ \hline\hline
1 & $g_2(z,a)=3(z-1)(z-a^2)^3$\\
& $g_3(z,a)=(z-1)(z-a^2)^4(z+a)$\\ \hline\hline
2 & $g_2(z,a)=12z^2(z^2+az+1)$\\
& $g_3(z,a)=4z^3(2z^3+3az^2+3az+2)$\\ \hline\hline
3 & $g_2(z,a)=12z^2(z^2+2az+1)$\\
& $g_3(z,a)=4z^3(2z^3+3(a^2+1)z^2+6az+2)$\\ \hline\hline
4 & $g_2(z,a)=3z^3(z+a)$\\
& $g_3(z,a)=z^5(z+1)$\\ \hline\hline
5 & $g_2(z,a)=3z^3(z+2a)$\\
& $g_3(z,a)=z^4(z^2+3az+1)$\\\hline\hline
\end{tabular}
\end{center}

\caption{The Herfurtner  list of "deformable" elliptic surfaces with four singular fibers}
\label{table6}
\end{table}
\begin{eqnarray*}
g_2=g_2(z,a)&=&3z^3(z+a)\\
g_3=g_3(z,a)&=&z^5(z+1) .
\end{eqnarray*}
We have
\begin{eqnarray*}
\Delta &=&\Delta (g_2,g_3)=27z^9((3a-2)z^2+(3a^2-1)z+a^3)\\
\delta &=&\delta (z,a)=-3z^7((3a-2)z+a).
\end{eqnarray*}
The Picard-Fuchs equation (\ref{pfeq}) takes the form

$$
144z^2((3a-2)z+a)((3a-2)z^2+(3a^2-1)z+a^3)^2 \, \eta_1''
$$
$$
+ 144z((3a-2)z^2+(3a^2-1)z+a^3)(3(3a-2)^2z^3
$$
$$
+2(3a-2)(3a-1)(a+1)z^2+a(3a^3+7a^2-3)z+2a^4)\, \eta_1'
$$
$$
+[135(3a-2)^3z^5+(3a-2)^2(468a^2+267a-164)z^4
$$
$$
+2(3a-2)(189a^4+522a^3-48a^2-208a+10)z^3
$$
$$
-2a(270a^5-1269a^4+252a^3+460a^2-70)z^2 $$ $$ -a^4(243a^3-666a^2+176)z+27a^7]\, \eta_1 = 0
$$ and  has four
regular singular points at $\infty$ and the roots of $(3a-2)z^2+(3a^2-1)z+a^3$ (the roots of $\Delta$), as well
one apparent singularity at the  root of $(3a-2)z+a $ (which is a root of $\delta$). Re-normalizing the singular
points to $0,1,t,\infty$ we get
\begin{equation}
 \lambda =\frac{a^2-a+1}{a^2(2-a)}, t=\frac{2a-1}{a^3(2-a)}, a\in \C . \label{4A}
\end{equation}
The parameter $a$ defines an algebraic isomonodromic deformation of the Picard-Fuchs equation (\ref{pfeq}) with
Riemann schema
$$
\begin{array}({ccccc})
0 & 1 & t & \lambda & \infty \\
-\frac{3}{4} & \frac{1}{4} & \frac{1}{4} & 0 & \frac{5}{4} \\
-\frac{1}{4} & -\frac{1}{4} & -\frac{1}{4} & 2 & \frac{3}{4}
\end{array} .
$$
Therefore the algebraic function $ \lambda =\lambda(t)$ determined implicitly by (\ref{4A}) is an algebraic
solution of $ \mathbf{PVI_\alpha}$ equation with
\begin{eqnarray*}
\alpha_0 =\frac{1}{8},\;\;\alpha_1 =\frac{1}{8}, \;\;\alpha_2 = \frac{1}{8},\;\;\alpha_3 =\frac{1}{8} .
\end{eqnarray*}
(see \cite{8} for details). Eliminating $a$ from  (\ref{4A}) we get
$$
\lambda^4-2t\lambda^3-2\lambda^3+6t\lambda^2-2t^2\lambda-2t\lambda+t^3-t^2+t =0
$$
which is an  equation for the solution 4A with $\alpha=(1/8,1/8,1/8,1/8)$. Together with the Doran's point
$\alpha= (0,1/8,0,0)$, this implies that the solution $\lambda(t)$ satisfies also the implicit equation
(\ref{1})
$$1 + \frac{t - 1}{(\lambda - 1)^2}
  -  \frac{t(t - 1)}{(\lambda - t)^2} = 0
$$
corresponding to the affine line through $(1/8,1/8,1/8,1/8)$ and $(0,1/8,0,0)$ described in Table \ref{table1},
4A. The solution (\ref{4A}) with $\alpha = (1/8,1/8,1/8,1/8)$ was found by Hitchin \cite[section 6.1]{7},
\cite[(34)]{h7}.
\begin{Rq}
If we repeat the same computation, but making use of the Picard-Fuchs system (\ref{pf1}) then of course we
obtain the same algebraic solution but with $\alpha= (1/8,0,0,0)$. This value has been erroneously computed by
Doran \cite{doran} to be  $(1/18,0,0,0)$. This led him to the wrong conclusion that the solution $4C$ is
equivalent by an Okamoto  transformation to the "cubic" solution $B_3$ of Dubrovin-Mazzocco \cite[p.140]{dub00}
with $\alpha= (25/18,0,0,0)$, see \cite[Remark 7]{doran}. As the Okamoto transformations of $
\mathbf{PVI_\alpha}$ act within the ring $\Z[1/2,\sqrt{2\alpha_1},
\sqrt{2\alpha_1},\sqrt{2\alpha_1},\sqrt{2\alpha_1}]$, then no solution of $
\mathbf{PVI_{\textit{(25/18,0,0,0)}}}$ is equivalent to a solution of $\mathbf{PVI_{\textit{(1/8,0,0,0)}}}$. M.
Mazzocco kindly informed us for a missprint in the formula for the $B_3$-solution, \cite[p.140]{dub00}. The
corrected formula is reproduced in \cite[formula (3.1)]{doran}.
\end{Rq}

\begin{table}

\begin{center}
\begin{tabular}{||c|c|c|c||}
\hline\hline
Stabilizer of  & Name of & $\textbf{PVI}_{\alpha}$ equation & Stabilizer of \\
the solution & the solution &  & $\textbf{PVI}_{\alpha}$ equation \\ \hline\hline
$D_4$ & $2A$ & ($0,0,\frac{1}{18},\frac{1}{18}$) & $S_2 \times S_2$ \\
 &  & ($\frac{1}{18},\frac{1}{18},0,0$) & \\
\cline{2-3}
 & $2B$ & ($\frac{1}{18},0,\frac{1}{18},0$) & \\
 &  & ($0,\frac{1}{18},0,\frac{1}{18}$) & \\
\cline{2-3}
 & $2C$ & ($\frac{1}{18},0,0,\frac{1}{18}$) & \\
 & & ($0,\frac{1}{18},\frac{1}{18},0$) & \\ \hline\hline
$D_4$ & $2A$ & ($0,0,0,0$) & $S_4$\\
 & $2B$ & & \\
 & $2C$ & & \\ \hline\hline
$S_3$ & $3A$ & ($0,0,0,0$) & $S_4$\\
 & $3B$ & & \\
 & $3C$ & & \\
 & $3D$ & & \\ \hline\hline
$S_3$ & $4A$ & ($0,\frac{1}{8},0,0$) & $S_3$\\
 & $4B$ & ($0,0,\frac{1}{8},0$) & \\
 & $4C$ & ($\frac{1}{8},0,0,0$) & \\
 & $4D$ & ($0,0,0,\frac{1}{8}$) & \\ \hline\hline
$S_2$ & $5A$ & ($0,\frac{1}{18},0,0$) & $S_3$\\
 & $5B$ & & \\
 & $5C$ & & \\
\cline{2-3}
 & $5D$ & ($0,0,\frac{1}{18},0$) & \\
 & $5E$ & & \\
 & $5F$ & & \\
\cline{2-3}
 & $5G$ & ($\frac{1}{18},0,0,0$) & \\
 & $5H$ & & \\
 & $5I$ & & \\
\cline{2-3}
 & $5J$ & ($0,0,0,\frac{1}{18}$) & \\
 & $5K$ & & \\
 & $5L$ & & \\ \hline\hline
\end{tabular}
\end{center}

\caption{Solutions $(\lambda(t),\alpha)$ of $ \mathbf{PVI_\alpha}$ equations related to the Picard-Fuchs system
(\ref{pf1}). }
 \label{table4}
\end{table}

\begin{table}

\begin{center}
\begin{tabular}{||c|c|c|c||}
\hline\hline
Stabilizer of  & Name of & $\textbf{PVI}_{\alpha}$ equation & Stabilizer of \\
the solution & the solution &  & $\textbf{PVI}_{\alpha}$ equation \\ \hline\hline
$D_4$ & $2A$ & ($\frac{1}{8},\frac{1}{8},\frac{1}{18},\frac{1}{18}$) & $S_2 \times S_2$\\
 &  & ($\frac{1}{18},\frac{1}{18},\frac{1}{8},\frac{1}{8}$) & \\
\cline{2-3}
 & $2B$ & ($\frac{1}{18},\frac{1}{8},\frac{1}{18},\frac{1}{8}$) & \\
 &  & ($\frac{1}{8},\frac{1}{18},\frac{1}{8},\frac{1}{18}$) & \\
\cline{2-3}
 & $2C$ & ($\frac{1}{18},\frac{1}{8},\frac{1}{8},\frac{1}{18}$) & \\
 & & ($\frac{1}{8},\frac{1}{18},\frac{1}{18},\frac{1}{8}$) & \\ \hline\hline
$D_4$ & $2A$ & ($\frac{1}{8},\frac{1}{8},\frac{1}{2},\frac{1}{2}$) & $S_2 \times S_2$ \\
 &  & ($\frac{1}{2},\frac{1}{2},\frac{1}{8},\frac{1}{8}$) & \\
\cline{2-3}
 & $2B$ & ($\frac{1}{2},\frac{1}{8},\frac{1}{2},\frac{1}{8}$) & \\
 &  & ($\frac{1}{8},\frac{1}{2},\frac{1}{8},\frac{1}{2}$) & \\
\cline{2-3}
 & $2C$ & ($\frac{1}{2},\frac{1}{8},\frac{1}{8},\frac{1}{2}$) & \\
 & & ($\frac{1}{8},\frac{1}{2},\frac{1}{2},\frac{1}{8}$) & \\ \hline\hline
$S_3$ & $3A$ & ($\frac{1}{8},\frac{9}{8},\frac{1}{8},\frac{1}{8}$) & $S_3$\\
 & $3B$ & ($\frac{9}{8},\frac{1}{8},\frac{1}{8},\frac{1}{8}$) & \\
 & $3C$ & ($\frac{1}{8},\frac{1}{8},\frac{9}{8},\frac{1}{8}$) & \\
 & $3D$ & ($\frac{1}{8},\frac{1}{8},\frac{1}{8},\frac{9}{8}$) & \\ \hline\hline
$S_3$ & $4A$ & ($\frac{1}{8},\frac{1}{8},\frac{1}{8},\frac{1}{8}$) & $S_4$\\
 & $4B$ & & \\
 & $4C$ & & \\
 & $4D$ & & \\ \hline\hline
$S_2$ & $5A$ & ($\frac{1}{2},\frac{1}{18},\frac{1}{8},\frac{1}{8}$) & $S_3$\\
 & $5B$ & ($\frac{1}{8},\frac{1}{18},\frac{1}{2},\frac{1}{8}$) & \\
 & $5C$ & ($\frac{1}{8},\frac{1}{18},\frac{1}{8},\frac{1}{2}$) & \\
 & $5D$ & ($\frac{1}{2},\frac{1}{8},\frac{1}{18},\frac{1}{8}$) & \\
 & $5E$ & ($\frac{1}{8},\frac{1}{2},\frac{1}{18},\frac{1}{8}$) & \\
 & $5F$ & ($\frac{1}{8},\frac{1}{8},\frac{1}{18},\frac{1}{2}$) & \\
 & $5G$ & ($\frac{1}{18},\frac{1}{8},\frac{1}{2},\frac{1}{8}$) & \\
 & $5H$ & ($\frac{1}{18},\frac{1}{2},\frac{1}{8},\frac{1}{8}$) & \\
 & $5I$ & ($\frac{1}{18},\frac{1}{8},\frac{1}{8},\frac{1}{2}$) & \\
 & $5J$ & ($\frac{1}{8},\frac{1}{8},\frac{1}{2},\frac{1}{18}$) & \\
 & $5K$ & ($\frac{1}{8},\frac{1}{2},\frac{1}{8},\frac{1}{18}$) & \\
 & $5L$ & ($\frac{1}{2},\frac{1}{8},\frac{1}{8},\frac{1}{18}$) & \\ \hline\hline
\end{tabular}
\end{center}

\caption{Solutions $(\lambda(t),\alpha)$ of $ \mathbf{PVI_\alpha}$ equations related to the Picard-Fuchs system
(\ref{pf2}). } \label{table5}
\end{table}

\vspace{2ex} \noindent
%%%%%%%%%%%%%

\end{document}